\documentclass[a4paper, final,11pt]{article}

\usepackage {amsmath,amsthm,amsfonts,amssymb,multicol}
\usepackage{graphicx}
\usepackage{booktabs}
\usepackage{array}
\usepackage{subfigure}
\usepackage[font=scriptsize]{caption}
\usepackage{enumerate}
\usepackage{url}
\usepackage[nodayofweek]{datetime}
\usepackage[margin=3cm]{geometry}
\usepackage[english]{babel}
\usepackage{verbatim} 
\usepackage{amsthm} 
\usepackage{enumitem}
\usepackage{enumerate}
\usepackage{bbm} 
\usepackage[normalem]{ulem} 
\usepackage{bm}
\usepackage[affil-it]{authblk}

\newcommand{\be}{\begin}
\newcommand{\e}{\end}
\newcommand{\beq}{\begin{equation}}
\newcommand{\eeq}{\end{equation}}
\newcommand{\beqs}{\begin{equation*}}
\newcommand{\eeqs}{\end{equation*}}
\newcommand{\bal}{\begin{align}}
\newcommand{\eal}{\end{align}}
\newcommand{\bals}{\begin{align*}}
\newcommand{\eals}{\end{align*}}

\renewcommand{\l}{\left}
\renewcommand{\r}{\right}

\newcommand{\set}[1]{\mathbb{#1}}

\newcommand{\curly}[1]{\mathcal{#1}}

\newcommand{\eps}{\epsilon}

\newcommand{\al}{\alpha}
\newcommand{\de}{\delta}

\newcommand{\tvector}[2]{\left(\be{array}{c}#1\\#2\e{array}\right)}

\renewcommand{\it}{\infty}

\theoremstyle{definition}

\theoremstyle{remark}

\def\dotuline{\bgroup
  \ifdim\ULdepth=\maxdimen  
   \settodepth\ULdepth{(j}\advance\ULdepth.4pt\fi
  \markoverwith{\begingroup
  \advance\ULdepth0.08ex
  \lower\ULdepth\hbox{\kern.15em .\kern.1em}%
  \endgroup}\ULon}

\def\dashuline{\bgroup
  \ifdim\ULdepth=\maxdimen  
   \settodepth\ULdepth{(j}\advance\ULdepth.4pt\fi
  \markoverwith{\kern.15em
  \vtop{\kern\ULdepth \hrule width .3em}%
  \kern.15em}\ULon}

\begin{document}

\title{New counterexamples for sums-differences}
\author{Marius Lemm\thanks{mlemm@caltech.edu}}
\affil{\small{Department of Mathematics, Caltech, Pasadena, CA 91125}}
\date{}

\maketitle

\abstract{We present new counterexamples, which provide stronger limitations to sums-differences statements than were previously known. The main idea is to consider non-uniform probability measures.}

\section{Introduction}

\paragraph{The sums-differences problem}
For $r\in \set{Q}\cup\{\it\}$, we define the maps $\pi_r:\set{R}^2\rightarrow \set{R}$ by
\beqs
    \pi_r(a,b)= a+rb,
\eeqs
with the convention that $a+\it b=b$.

Let $r_1,\ldots,r_n\in \set{Q}\cup\{\it\}\setminus\{-1\}$ and $1<\al\leq 2$. We say that the statement $SD(r_1,\ldots,r_n;\al)$ holds, if for any number $N$ and any finite $G\subset \set{R}^2$ such that $\pi_{-1}$ is injective on $G$ and $|\pi_{r_j}(G)|\leq N$ for $j=1,\ldots,n$, it must be that $|G|<N^\al$. Here, $|A|$ denotes the cardinality of a set $A$.

We say that the statement $SD(\al)$ holds, if for any $\eps>0$, there exists $r_1,\ldots,r_n\in \set{Q}\cup\{\it\}\setminus\{-1\}$ such that $SD(r_1,\ldots,r_n;\al+\eps)$.\\

The sums-differences problem was introduced by Bourgain in \cite{Bourgain99} as a combinatorial tool to prove lower bounds on the Hausdorff dimension of Kakeya sets in $\set{R}^d$. Already the trivial statement $SD(2)$ implies the a priori non-trivial fact that the Hausdorff dimension of any Kakeya set is at least $(d+1)/2$. A proof of $SD(1)$ would imply the Kakeya conjecture, i.e.\ that any Kakeya set has Hausdorff dimension $d$. The goal of the sums-differences approach to the Kakeya problem is then to prove $SD(\al)$ statements with $\al$ as small as possible.

In \cite{Bourgain99}, Bourgain showed the statement $SD(0,1,\it;25/13)$. In \cite{KatzTao99}, Katz and Tao proved $SD(0,1,\it;11/6)$ and $SD(0,1,2,\it;7/4)$. In \cite{KatzTao02}, they improved to $SD(\al)$ with $\al\approx 1.67513$ using iteration. Well over ten years later, these are still the best known values of $\al$.\\

Our results will go in the opposite direction. We say $\neg SD(r_1,\ldots,r_n;\al)$, if there exists a counterexample, i.e.\ an explicit choice of $G$ such that $\pi_{-1}$ injective on $G$ and such that
\beq
\label{eq:cardal}
    \al = \frac{\log |G|}{\max_j\{\log|\pi_{r_j}(G)|\}}.
\eeq
In this paper, we present new counterexamples that prove $\neg SD(0,1,\it;\al)$ as well as $\neg SD(0,1,2,\it;\al)$ for considerably larger values of $\al$ than previously known.

\paragraph{Entropy formulation of sums-differences }
We will state our counterexamples in an equivalent formulation of the problem, in which the logarithm of cardinality in \eqref{eq:cardal} is replaced by entropy, see \eqref{eq:entropyal}. The entropy viewpoint is instrumental to our construction of counterexamples. It is due to Ruzsa, but was not published by him and is therefore not as well-known as it could be. We discuss it here in some detail.

Let $I$ be a finite set. We denote by $\curly{M}(I)$ the set of probability measures on $I$. Given $P\in \curly{M}(I)$, we define its \emph{entropy} $H(P)$ by
\beqs
    H(P) = - \sum_{i\in\set{I}} p_i \log p_i.
\eeqs
Let $J$ be another finite set. Given a map $f:I\rightarrow J$, the \emph{push-forward measure} $fP$ is defined by
\beqs
fP(j)=P(f^{-1}(j))
\eeqs
 for all $j\in J$.

\be{prop}[Entropy formulation, Ruzsa]
\label{prop:entropy}
Let $r_1,\ldots,r_n\in \set{Q}\cup\{\it\}\setminus\{-1\}$.\\ The following are equivalent:
\be{enumerate}[label=(\roman*)]
 \item $\neg SD(r_1,\ldots,r_n;\al)$
 \item there exists a finite set $G\subset\set{R}^2$ such that $\pi_{-1}$ is injective on $G$ and $P\in \curly{M}(G)$ such that
\beq
\label{eq:entropyal}
    \frac{H(P)}{\max_j\{H(\pi_{r_j}P)\}}\geq \al.
\eeq
\e{enumerate}
\e{prop}

\be{proof}

\dashuline{(ii)$\Leftarrow$ (i):} We will use without proof the fact that the sums-differences problem is independent of the choice of underlying vector space, i.e.\ instead of $G\subset\set{R}^2$ one can equivalently consider $G\subset V^2$ for any vector space $V$. The reason for this is that one can re-formulate sums-differences as a purely graph-theoretical problem, see \cite{Katz06}.

Suppose we are given a finite set $G\subset\set{R}^2$ such that $\pi_{-1}$ is injective on $G$ and $P\in \curly{M}(G)$ such that \eqref{eq:entropyal} holds. For every $\eps>0$, we will construct an explicit $G'\subset\l(\set{R}^{M}\r)^2$ for some large $M$ such that $\pi_{-1}$ is injective on $G'$ and such that
\beq
\label{eq:claim1}
    \al':= \frac{\log|G'|}{\max_j\{\log|\pi_{r_j}(G')|\}} > \al - \eps.
\eeq
Taking $\eps\rightarrow 0$ then proves $\neg SD(r_1,\ldots,r_n;\al)$.

The basic idea of our construction is that we can approximate any $P$ by a multinomial distribution that arises from appropriately binning uniform measure on a large underlying set.

Let $\de>0$. For all $g\in G$, we can find a rational number $q_g$ such that $|q_g-P(g)|<\de$. We may arrange that $\sum_g q_g=1$. Let $M$ denote the largest denominator appearing in the collection $\{q_g\}_g$ after maximal reduction. Let $k_g$ denote the positive integer $q_g M$. We define $G'\subset \l(\set{R}^{M}\r)^2$ by
\begin{align*}
    G'=&\left\{ \tvector{\l(g_1(1),\ldots,g_M(1)\r)}{\l(g_1(2),\ldots,g_M(2)\r)}\in \l(\set{R}^{M}\r)^2 :\right. \\
      &\;\left.\vphantom{\tvector{\l(g_1(1),\ldots,g_M(1)\r)}{\l(g_1(2),\ldots,g_M(2)\r)}} g_i=\l(g_i(1),g_i(2)\r)\in G, \forall i \wedge |\{i: g_i = g\}|=k_g  , \forall g\in G\right\}.
\end{align*}
and observe that
\beqs
    |G'|= \frac{M!}{\prod_{g\in G} k_g!}.
\eeqs
It is easy to check that
\begin{align*}
    &\pi_{r_j}(G')\\
     &= \l\{ \l(\nu_1,\ldots,\nu_M\r)\in \set{R}^M : \nu_i\in \pi_{r_j}(G), \forall i \wedge |\{i: \nu_i = \nu\}|=\sum_{g\in \pi_{r_j}^{-1}(\nu)} k_g , \forall \nu\in \pi_{r_j}(G)\r\},
\end{align*}
and so
\beqs
    |\pi_{r_j}(G')| = \frac{M!}{\prod_{\nu\in \pi_{r_j}(G)} \l(\sum_{g\in \pi_{r_j}^{-1}(\nu)} k_g \r)!}.
\eeqs
In particular, $|\pi_{-1}(G)|=|G|$ implies $|\pi_{-1}(G')|=|G'|$. We now have everything to compute
\beqs
    \al' = \frac{\log|G'|}{\max_j\{\log|\pi_{r_j}(G')|\}} = \frac{\log(M!)-\sum_{g\in G} \log(k_g!)}{\max_j\l\{\log(M!)-\sum_{\nu\in \pi_{r_j}(G)} \log\l(\l(\sum_{g\in \pi_{r_j}^{-1}(\nu)} k_g \r)!\r)\r\}}.
\eeqs
The key fact that (asymptotically) relates combinatorics to entropy is Stirling's formula: $\log(N!)=N \log(N/e)+O(N)$. Introducing $\psi(x)=-x\log x$, we obtain
\beqs
    \al'=\frac{-\sum_{g\in G} \psi\l(\frac{k_g}{M}\r)}{\max_j\l\{\sum_{\nu\in \pi_{r_j}(G)} \psi\l(\sum_{g\in \pi_{r_j}^{-1}(\nu)}\frac{k_g}{M}\r)\r\}} + o(1).
\eeqs
Since $k_g/M=q_g$ can be made arbitrarily close to $P(g)$, \eqref{eq:claim1} follows.

\dashuline{ (i)$\Rightarrow$ (ii):} This direction is easy: Take $P$ to be uniform measure on the finite set $G\subset\set{R}^2$ provided by $\neg SD(r_1,\ldots,r_n;\al)$. Then, $H(P)=\log(|G|)$. By Jensen's inequality, $H(\pi_{r_j}P)\leq \log(|\pi_{r_j}|)$ and so
\beqs
    \frac{H(P)}{\max_j\{H(\pi_{r_j}P)\}} \geq \frac{\log(|G|)}{\max_j\{\log(|\pi_{r_j}|)\}} =\al.
\eeqs
\e{proof}

\paragraph{Ruzsa's counterexample}
The following classical construction is due to Ruzsa. It yields $\neg SD(0,1,\it;\frac{\log(27)}{\log(27/4)})$ and was the best known counterexample in that case so far. Let
\beqs
    G=
    \l\{\begin{array}{c}
        (0,1)\\
        (1,0)\\
        (1,1)
    \end{array}\r\}
\eeqs
and note that $|\pi_{-1}(G)|=|G|=3$. Let $P$ be uniform probability measure on $G$, i.e.\ $P$ assigns probability $1/3$ to each element of $G$. Then, an easy computation shows
\beqs
    \frac{H(P)}{\max_j\{H(\pi_{r_j}P)\}} = \frac{\log(27)}{\log(27/4)}.
\eeqs
By Proposition \ref{prop:entropy}, this implies $\neg SD(0,1,\it;\frac{\log(27)}{\log(27/4)})$.

\paragraph{Acknowledgement} The author wishes to thank Nets Katz for encouragement and advice.

\section{Results}
To motivate our first result, let us compare Ruzsa's $\frac{\log(27)}{\log(27/4)}\approx 1.726$ with $11/6\approx 1.833$, the best known value of $\al$ for which $SD(0,1,\it;\al)$ is known to hold. There is a gap of size $\approx 0.1$ between these two values. Since improvements over Ruzsa's counterexample were elusive, it was believed that $SD(0,1,\it;\al)$ could hold for all $\al>\frac{\log(27)}{\log(27/4)}$.

\be{thm}
\label{thm1}
There exists $\al>1.77898$, such that $\neg SD(0,1,\it,\al)$.
\e{thm}

This manages to close about half of the $0.1$ gap. The value of $1.77898$ is obtained by numerical nonlinear maximiziation, which is notorious for getting stuck in local extrema. Thus, it is not clear that this value is best possible.

We also have
\be{thm}
\label{thm2}
There exists $\al>1.61226$, such that $\neg SD(0,1,2,\it,\al)$.
\e{thm}

The number $1.61226$ is to be compared with $7/4$, the best known value of $\al$ for which $SD(0,1,2,\it;\al)$ is known to hold.

\paragraph{The main idea}
To explain our approach, it is instructive to consider a simple modification of Ruzsa's construction that already yields $\neg SD(0,1,\it,\al)$ with an explicit $\al$ satisfying $\al\approx1.7726$. We will see that choosing a (particular) \emph{non-uniform} $P$ is what enables us to improve over Ruzsa's original counterexample, which featured uniform measure. Let
\beqs
    G=
    \l\{\begin{array}{c}
        (0,1)\\
        (1,0)\\
        (1,1)\\
        (2,0)
    \end{array}\r\}
\eeqs
and note that $|\pi_{-1}(G)|=|G|=4$. According to Proposition \ref{prop:entropy}, we can consider any probability measure $P=(p_1,p_2,p_3,p_4)$ on $G$.

Introducing the functions
\beq
\label{eq:psidefn}
    \psi(x)=-x\log x,\qquad \phi(x)=\psi(x)+\psi(1-x)
\eeq
the entropies can be written as
\begin{align*}
H(P)&=\psi(p_1)+\psi(p_2)+\psi(p_3)+\psi(1-p_1-p_2-p_3)\\
H(\pi_0P)&=\psi(p_2+p_3)+\psi(p_1)+\psi(1-p_1-p_2-p_3)\\
H(\pi_1P)&=\phi(p_1+p_2)\\
H(\pi_{\it}P)&=\phi(p_1+p_3).
\end{align*}
We aim to find the values of $p_1,p_2,p_3$ that maximize
\beqs
    \al(P)= \min_{j=0,1,\it} \frac{H(P)}{H(\pi_jP)}.
\eeqs
The minimum renders this non-differentiable and prevents it from being a calculus problem. From symmetry and convexity considerations, it is sensible to set
\beqs
H(\pi_0P)=H(\pi_1P)=H(\pi_\it P).
\eeqs
Elementary computation then yields
\beqs
     p_3=p_2,\quad p_2=1/2-p_1, \quad \log(2)=\psi(1-2p_1)+2\psi(p_1)
\eeqs
which uniquely determine $p_1,p_2,p_3$ and thus $\al$. The numerical values are $p_1\approx .1135$, $p_2=p_3\approx .3865$ and $\al\approx 1.772$.

\paragraph{Proofs}
We know give the counterexamples that prove Theorems \ref{thm1} and \ref{thm2}.

\be{proof}[Proof of Theorem \ref{thm1}]
We consider
\beqs
    G=
    \l\{\begin{array}{c}
        (0,1)\\
        (1,1)\\
        (1,0)\\
        (2,0)\\
        (2,-1)\\
        (3,-1)\\
        (3,-2)
    \end{array}\r\}.
\eeqs
Note that $|\pi_{-1}(G)|=|G|$. If we denote $P=(p_1,\ldots,p_7)$, we can set $p_7=p_1,\; p_6=p_2,\; p_5=p_3$
 to ensure that $H(\pi_0P)=H(\pi_\it P)$. By numerical maximization of
\beqs
    \frac{H(P)}{\max\{H(\pi_0P),H(\pi_1P)\}},
\eeqs
we find that the choice $p_1\approx .00024983, p_2\approx .028156, p_3\approx .22425 $ together with Proposition \ref{prop:entropy} imply Theorem \ref{thm1}.
\e{proof}

\be{rmk} It is noteworthy that, when one adds the points $(4,-2),(4,-3)$ to $G$ and sets $p_i=p_{9-i}$ for all $i$, one does not obtain a better value of $\al$, at least on the level of numerics.
\e{rmk}

\be{proof}[Proof of Theorem \ref{thm2}]
We consider
\beqs
    G=
    \l\{\begin{array}{c}
        (0,1)\\
        (1,0)\\
        (1,1)\\
        (2,0)\\
        (1,\frac{1}{2})
    \end{array}\r\}.
\eeqs
Note that $|\pi_{-1}(G)|=|G|$. Setting $p_1=p_2=p_3=p_4=: p$ ensures
\begin{align*}
H(\pi_1P)=H(\pi_\it P),\quad H(\pi_0P)=H(\pi_2 P).
\end{align*}
Furthermore, we suppose that $H(\pi_1P)=H(\pi_0P)$ or equivalently
\beqs
    2\psi(2p)+\psi(1-4p) =\psi(1-2p)+2\psi(p)
\eeqs
with $\psi$ as in \eqref{eq:psidefn}. This has a unique non-zero solution which satisfies $p\approx .21798$. By Proposition \ref{prop:entropy}, Theorem \ref{thm2} follows.
\e{proof}

\bibliographystyle{amsplain}
\bibliography{final}

\end{document}